\documentclass[a4paper,twoside,fullpage]{article}
\newtheorem{thm}{Theorem}[section]
\newtheorem{defn}{Definition}[section]
\def\halmos{\hfill\rule{6pt}{6pt}}
\newcommand{\kerR}{\mathop{\rm ker_{\ssR}}\nolimits}
\newcommand{\ssR}{\scriptscriptstyle R}
\newcommand{\kerL}{\mathop{\rm ker_{\ssL}}\nolimits}
\newcommand{\ssL}{\scriptscriptstyle L}

\usepackage{amssymb}
\makeatletter
\@addtoreset{equation}{section}
\def\thethm{\arabic{section}.\arabic{thm}\protect\@blinkpoint}
\def\theprop{\arabic{subsection}.\arabic{prop}\protect\@blinkpoint}
\def\thecor{\arabic{subsection}.\arabic{cor}\protect\@blinkpoint}
\def\thedefn{\arabic{section}.\arabic{defn}\protect\@blinkpoint}

\def\therem{\arabic{subsection}.\arabic{rem}\protect\@blinkpoint}
\def\thelemm{\arabic{subsection}.\arabic{lemm}\protect\@blinkpoint}
\def\theexmp{\arabic{subsection}.\arabic{exmp}\protect\@blinkpoint}
\def\thesection{\arabic{section}\protect\@blinkpoint}
\def\thesubsection{\arabic{subsection}\protect\@blinkpoint}
\def\thesubsubsection{\arabic{subsection}.\arabic{subsubsection}\protect\@blinkpoint}
\def\@blinkpoint{.}
\let\@blinkref=\ref
\def\ref#1{{\def\@blinkpoint{}\@blinkref{#1}}}
\let\@afterindentfalse\@afterindenttrue
\makeatother

\font\msbm=msbm10 scaled 1200
\font\msbmscript=msbm8
\font\msbmscriptscript=msbm6
\newfam\Bbbfam
\textfont\Bbbfam=\msbm
\scriptfont\Bbbfam=\msbmscript
\scriptscriptfont\Bbbfam=\msbmscriptscript
\def\Bbb#1{{\fam\Bbbfam#1}}

\begin{document}

\begin{center}
{\Large \bf Generalized inversion of Toeplitz-plus-Hankel
matrices}

V.\ M.\ Adukov, O.\ L.\ Ibryaeva

avm@susu.ac.ru, \ oli@susu.ac.ru

{\em South Ural State University, Chelyabinsk, Russia.}
\end{center}

A generalized inversion of block $T+H$ matrix is obtained for the 
first time. In a particular case when $T+H$ matrix is invertible, 
the method allows to obtain its inverse matrix without an additional
condition of invertibility of the corresponding $T-H$ matrix.

Keywords:  Toeplitz-plus-Hankel matrices, generalized inversion.

\section{Introduction}
\indent
In many applications, e.g. digital signal processing, discrete inverse
scattering, linear prediction etc.,
Toeplitz-plus-Hankel ($T+H$) matrices need to be inverted.
(For further applications see \cite{Bruckstein} and references therein).

Firstly the $T+H$ matrix inversion problem has been solved in
\cite{Merchant} where it was reduced to the inversion problem
of the block Toeplitz matrix (the so-called mosaic matrix).
The drawback of the method is that it does not work for any invertible
$T+H$ matrix since it requires also invertibility of the corresponding
$T-H$ matrix. This drawback appeared also in \cite{Ners}, \cite{Heinig-Rost-book}.
Later on the drawback was put out (see, e.g. \cite{Heinig-Rost-book}, \cite{Heinig-88}),
moreover, the inversion problem was solved for the block $T+H$ matrix
\cite{Shalom},\cite{Ad-Ib}.

The generalized inversion of Toeplitz-plus-Hankel matrices is of great interest
to, e.g., Pade-Chebyshev approximations. The generalized inversion for matrix $A$ is meant to be the matrix
$A^{\dagger}$ such that $AA^{\dagger}A=A.$
Our goal is to obtain the generalized inversion of block $T+H$ matrices
with help of the well-known method of reducing the block-diagonal matrix
formed from the $T+H$ and $T-H$ matrices to the mosaic matrix (as in the
work \cite{Merchant}).

We will need the generalized inversion for the block Toeplitz matrix which
has been already found in, e.g. \cite{Adukov-98}.
It is shown in the present paper that there is no need for $T-H$ matrix to be inverted: if the $T+H$ matrix
is invertible than the obtained generalized inverse matrix proves to be
its inverse matrix.

The paper is organized as follows. At first the basic definitions and the
main results of the paper \cite{Adukov-98} are given, then our main theorem
is formulated and proved. This theorem is demonstrated with an example in the
end of the paper.

This work was supported by Russian Foundation for Basic Research (RFFI),
grant N 04-04-96006.

\section{The basic definitions and notations}
\indent
We are going to find a generalized inversion for the block Toeplitz-plus-Hankel
matrix
\[
T+H=
\left(
\begin{array}{cccc}
a_{0} & a_{-1} & \ldots& a_{-m} \\
a_{1} & a_0 & \ldots &a_{-m+1}\\
\vdots &\vdots& \ddots &\vdots \\
a_{n} & a_{n-1}& \ldots & a_{n-m}
\end{array}
\right)
+
\left(
\begin{array}{cccc}
b_0 & b_{1}&\ldots& b_m \\
b_1 & b_{2}&\ldots &b_{m+1}\\
\vdots & \vdots& \ddots &\vdots \\
b_n & b_{n+1}&\ldots & b_{n+m}
\end{array}
\right),
\]
with $a_j,b_j \in \Bbb C^{p\times q}$.

We denote $a_{-m}^n(z)=a_{-m}z^{-m}+\ldots+a_0+\ldots+a_nz^n,$
$b_0^{n+m}(z)=b_0+b_1z+\ldots+b_{n+m}z^{n+m}$
and introduce an auxilary matrix function
\[A(z)=
\left(
\begin{array}{cc}
z^nb_0^{n+m}(z^{-1}) & z^{n-m}a_{-m}^n(z^{-1}) \\
a_{-m}^n(z) & z^{-m}b_0^{n+m}(z)
\end{array}
\right).
\]

Obviously,  $A(z)=\sum\limits_{j=-m}^n A_jz^j,$ with
$A_j \in \Bbb C^{2p \times 2q}$ and

\begin{equation}
\label{A}
A_j=
\left(
\begin{array}{cc}
b_{n-j} & a_{n-m-j}\\
a_j& b_{j+m}
\end{array}
\right).
\end{equation}

Thus, $A(z)$ is the generating function for the sequence of matrices
$A_{-m}, \ldots, A_0,$ $ \ldots, A_n$.

Later on we will need a generalized
inversion for the block Toeplitz matrix
\[
T_A=
\left(
\begin{array}{cccc}
A_0 & A_{-1} & \ldots & A_{-m}\\
A_{1} & A_0 & \ldots & A_{-m+1}\\
\vdots & \vdots &\ddots &\vdots \\
A_n& A_{n-1}& \ldots & A_{n-m}
\end{array}
\right).
\]

It has been already found in \cite{Adukov-98}. In order to use this result
we should introduce the notations of the essential indices and polynomials
of the sequence $A_{-m}, \ldots, A_0, \ldots, A_n$ .

We include the matrix $T_A\equiv T_0$ into the family of the
block Toeplitz matrices
\[T_k=
\left(
\begin{array}{cccc}
A_k & A_{k-1} & \ldots & A_{-m}\\
A_{k+1} & A_k & \ldots & A_{-m+1}\\
\vdots & \vdots &\ddots &\vdots \\
A_n& A_{n-1}& \ldots & A_{n-m-k}
\end{array}
\right), \ \ -m\leq k\leq n.
\]

The matrices $T_k$ are of the same structure and it is reasonable that they should
be examined together.

We are interested in right kernels of $T_k$.
For the sake of convenience let us pass from the spaces $\kerR T_k$ to the
isomorphic spaces ${\cal N}_k^{\ssR}$ of generating polynomials.

To do this we define the operator $\sigma_{\ssR}$ acting from
the space of rational matrix functions
$R(z)=\sum_{j=-n}^m r_jz^j, r_j \in \Bbb C^{2q\times l}$
to the space $\Bbb C^{2p\times l}$ according to
\[
\sigma_{\ssR}\left\{R(z)\right\}=\sum_{j=-n}^m A_{-j}r_j.
\]

By ${\cal N}_k^{\ssR}, k=-m,\ldots,n$, we denote the space of vector
polynomials $R(z)=\sum_{j=0}^{k+m} r_jz^j, r_j \in \Bbb C^{2q\times 1}$,
such that
\[
\sigma_{\ssR}\left\{z^{-i}R(z)\right\}=0, \ \ \ i=k,k+1,\ldots, n.
\]

${\cal N}_k^{\ssR}$ is evident to be isomorphic to $\kerR T_k.$

It is convenient to put ${\cal N}_{-m-1}^{\ssR}=0$ and denote
by ${\cal N}_{n+1}^{\ssR}$ the $2(n+m+2)q$-dimensional
space of all vector polynomials in $z$ with formal degree $n+m+1$.

Similarly, one may define spaces ${\cal N}_k^{\ssL}$ which are isomorphic to $\kerL T_k.$
We denote $ \ker_{L}A = \{y\ |\ yA=0\}.$

Let us put also $\alpha=\dim{\cal N}_{-m}^{\ssR}$ and $\omega=\dim{\cal N}_{n}^{\ssL}.$

We will say that the sequence $A_{-m}, \ldots, A_n$ is {\it left (right) regular}
if $\alpha=0\ \ (\omega=0).$  Otherwise, the sequence is not regular
and $\alpha (\omega)$ is its {\it left (right) defect}.
The sequence is called {\it regular} if $\alpha=\omega=0$.
It is evident that $\alpha < 2q, \ \ \omega< 2p$ for the nonzero sequence.

We will denote $d_k^{\ssR}=\dim{\cal N}_k^{\ssR}, \ \ $
$\Delta_k^{\ssR}=d_k^{\ssR}-d_{k-1}^{\ssR}, \ \ k=-m, \ldots, n+1$.

As it is proved in \cite{Adukov-98}, for any sequence $A_{-m}, \ldots, A_n$
the following inequalities hold:
\[
\alpha=\Delta_{-m}^{\ssR}\leq \Delta_{-m+1}^{\ssR}\leq \ldots\leq
\Delta_n^{\ssR}\leq \Delta_{n+1}^{\ssR}=2(p+q)-\omega.
\]
It means that there are $2(p+q)-\alpha-\omega$ integers $\mu_{\alpha}
\leq \mu_{\alpha+1}\leq \ldots \leq \mu_{2(p+q)-\alpha-\omega}$, satisfying
equations
\begin{equation}
\label{table}
\begin{array}{ccccccl}
\Delta^{\ssR}_{-m}&=&\ldots&=&\Delta^{\ssR}_{\mu_{\alpha+1}}&=
&\alpha, \\
& &\cdots& & & &   \\
\Delta^{\ssR}_{\mu_i+1}&=&\ldots&=&\Delta^{\ssR}_{\mu_{i+1}}&=&i, \\
& &\cdots& & & &  \\
\Delta^{\ssR}_{\mu_{2(p+q)-\omega}+1}&=&\ldots&=&\Delta^{\ssR}_{n+1}&=
&2(p+q)-\omega.
\end{array}
\end{equation}
If the $i$th row in (\ref{table}) is absent, we assume $\mu_i=\mu_{i+1}.$
Let us put also $\mu_1=\ldots=\mu_{\alpha}=-m-1$ if $\alpha \ne 0$
and $\mu_{2(p+q)-\omega+1}=\ldots=\mu_{2(p+q)}$ if $\omega \ne 0.$

Thus there is a set of $2(p+q)$ integers, satisfying (\ref{table}), for
any sequence $A_{-m},\ldots,A_n$.
We will call these integers as {\em indices} of the sequence.

Now we will define the right essential polynomials of the
sequence.

It follows from the definition of ${\cal N}_k^{\ssR}$ that
 ${\cal N}_k^{\ssR}$ and $z{\cal N}_k^{\ssR}$ are the subspaces
of ${\cal N}_{k+1}^{\ssR}$, $k=-m-1,\ldots,n$,
moreover, ${\cal N}_k^{\ssR}\bigcap z{\cal N}_k^{\ssR}={\cal N}_{k-1}^{\ssR}$.
Then
\[
{\cal N}_{k+1}^{\ssR}=\left({\cal N}_k^{\ssR}+z{\cal N}_k^{\ssR}\right)\oplus
{\cal H}_{k+1}^{\ssR},
\]
where ${\cal H}_{k+1}^{\ssR}$ is the complement of
${\cal N}_k^{\ssR}+z{\cal N}_k^{\ssR}$ to the whole ${\cal N}_{k+1}^{\ssR}.$

Obviously,  $\dim {\cal H}_{k+1}^{\ssR}=\Delta_{k+1}^{\ssR}-\Delta_{k}^{\ssR}.$
Hence $\dim {\cal H}_{k+1}^{\ssR}\ne 0$ iff $k=\mu_i.$
In this case $\dim {\cal H}_{k+1}^{\ssR}$ is equal to the multiplicity
$k_i$ of the index $\mu_i.$

\begin{defn}
If $\alpha\ne 0$ then any column polynomials
$R_1(z), \ldots,$ $R_{\alpha}(z)$ forming the basis of ${\cal N}_{-m}^{\ssR}$
will be called right essential polynomials of the sequence
$A_{-m},\ldots,$ $A_0,$ $\ldots,A_n$. They correspond to the index $\mu_1=-m-1$ with
the multiplicity $\alpha$.

Any vector polynomials $R_j(z), \ldots, R_{j+k_j-1}(z)$ forming the basis
for ${\cal H}_{\mu_j+1}^{\ssR}$ will be called right essential polynomials of the sequence
$A_{-m},\ldots,$ $A_0,$ $\ldots,A_n$. They correspond to the index $\mu_j$
with the multiplicity $k_j$, $\alpha+1 \leq j\leq 2(p+q)-\omega.$
\end{defn}

Similarly, one may define the left essential polynomials.

There are
$2(p+q)-\omega$ right and $2(p+q)-\alpha$ left essential polynomials
of the sequence $A_{-m},\ldots,A_n$.

If $\alpha \ne 0$ or $\omega \ne 0$ then there is a lack of essential polynomials.
But actually we can always complement the number of right (when $p\leq q$)
or left  (when $p\geq q$) essential polynomials to $2(p+q).$
(The complement procedure was described in \cite{Adukov-98}).

Henceforth,  for definiteness sake, we will suppose that we have got the full set of
$2(p+q)$ right essential polinomials, i.e. either $\omega=0$ or $p\leq
q$.

The set of the left essential polynomials could always be recovered with
the help of the so-called conformation procedure of the right and left essential
polynomials.
Let us describe how for the given set of the right essential polynomials
$R_1(z),$ $\ldots,R_{2(p+q)}(z)$, $R_j(z)\in {\Bbb C}^{2q\times 1}[z]$
one can construct the conforming left essential polynomials $L_1(z),\ldots,L_{2(p+q)}(z)$,
 $L_j(z)\in {\Bbb C}^{1\times 2p}[z]$.

We introduce the matrix of the right essential polynomials
\[
{\cal R}(z)=\left(R_1(z) \ldots R_{2(p+q)}(z)\right)
\]
and find the matrix polynomial $\alpha_{-}(z)$ from the next decomposition
\[
A(z){\cal R}(z)=\alpha_{-}(z)d(z)-z^{n+1}\beta_{+}(z).
\]
Here $d(z)={\rm diag} \left[z^{\mu_1},\ldots,z^{\mu_{2(p+q)}}\right],$
$\beta_{+}(z)$ is the matrix polynomial in $z$, $\alpha_{-}(z)$ is
the matrix polynomial in $z^{-1}$. Both of them have sizes $2p\times 2(p+q)$.

Let $U_{-}(z)$ be the matrix polynomial in $z^{-1}$ such that:
\[
U_{-}(z)=
\left(
\begin{array}{c}
{\cal R}_{-}(z) \\
\alpha_{-}(z)
\end{array}
\right),
\]
with ${\cal R}_{-}(z)=z^{-m-1}{\cal R}(z)d^{-1}(z).$

The matrix polynomial $U_{-}(z)$ is shown in \cite{Adukov-98} to be unimodular, i.e.
its determinant is equal to a constant.

We pick the  $2(p+q)\times 2p$ block ${\cal L}(z)$ out $U^{-1}_{-}(z)$:
\[
U^{-1}_{-}(z)=\left(
\begin{array}{cc}
* & {\cal L}(z)
\end{array}
\right).
\]

The matrix polynomial
\[
{\cal L}(z)=
\left(
\begin{array}{c}
L_1(z) \\
\vdots \\
L_{2(p+q)}(z)
\end{array}
\right)
\]
turns out to be the matrix of the conforming left essential polynomials.

The case when $\alpha =0$ or $p\geq q$ may be considered in a similar manner
with help of the left essential polynomials.

Now we may present the formula (5.13) from \cite{Adukov-98}
for the generalized inverse of $T_A$:
\begin{equation}
\label{513}
T_A^{\dagger}=
\left(
\begin{array}{ccc}
{\cal R}_0 & \ldots & 0 \\
\vdots &\ddots & \vdots \\
{\cal R}_m & \ldots & {\cal R}_0
\end{array}
\right)\Pi
\left(
\begin{array}{ccc}
{\cal L}_0 & \ldots &  {\cal L}_{-n} \\
\vdots &\ddots & \vdots \\
0 & \ldots & {\cal L}_{0}
\end{array}
\right).
\end{equation}

Here ${\cal R}_j \in {\Bbb C}^{2q\times 2(p+q)},$ \ \
${\cal L}_j \in {\Bbb C}^{2(p+q)\times 2p}$ are the coefficients of the
matrix polynomials ${\cal R}(z)$, ${\cal L}(z)$, respectively, and
$R_j(z), L_j(z)$ are the conforming right and left essential polynomials
of the sequence $A_{-m}, \ldots, A_0,$ $ \ldots, A_n$.
The generalized inversion for matrix $A$ is meant to be the matrix $A^{\dagger}$ such
that $AA^{\dagger}A=A.$

The matrix $\Pi$ is constructed as follows.
Let $\lambda_1,\ldots,\lambda_r$ be the distinct essential indices
of the sequence $A_{-m}, \ldots, A_0, \ldots, A_n$ and let
$\nu_1,\ldots,\nu_r$ be their multiplicities ($\nu_1+\ldots+\nu_r=2(p+q)$).

Then
\[
\Pi=
\left(
\begin{array}{cccc}
\Pi_0 & \Pi_{-1}&\ldots& \Pi_{-n} \\
\Pi_1 & \Pi_{0}&\ldots &\Pi_{-n+1}\\
\vdots & \vdots& \ddots &\vdots \\
\Pi_m & \Pi_{m-1}&\ldots & \Pi_{m-n}
\end{array}
\right).
\]
Here $\Pi_k=0$ for
$-n\leq k\leq m,\ k\neq -\lambda_1,\ldots,-\lambda_r$,
 $\Pi_{-\lambda_j}=\|\varepsilon_i^j\delta_{ik} \|_{i,k=1}^{2(p+q)}$,
\[
\varepsilon_i^j=
\left\{
\begin{array}{cc}
1, & i=\nu_1+\cdots+\nu_{j-1}+1,\ldots,\nu_1+\cdots+\nu_j, \\
0, & {\rm otherwise}.
\end{array}
\right.
\]

The following partition of the right essential polynomials $R_j(z)$  will be useful
for the generalized inversion of the $T+H$ matrix:
\[
R_j(z)=
\left(
\begin{array}{c}
R_j^1(z) \\
\\
R_j^2(z)
\end{array}
\right).
\]
Here  $R_j^{1,2}\in {\Bbb C}^{q\times 1}[z].$ In similar way we partition
the left essential polynomials:
\[
L_j(z)=
\left(
\begin{array}{cc}
L_j^1(z) & L_j^2(z)
\end{array}
\right),
\]
with  $L_j^{1,2}\in {\Bbb C}^{1\times p}[z].$

Then the matrix of these essential polynomials may be represented as:
\begin{equation}
\label{div}
{\cal R}(z)=
\left(
\begin{array}{c}
{\cal R}^1(z) \\
{\cal R}^2(z)
\end{array}
\right), \ \ \
{\cal L}(z)=
\left(
\begin{array}{cc}
{\cal L}^1(z) & {\cal L}^2(z)
\end{array}
\right),
\end{equation}
with ${\cal R}^{1,2}(z) \in {\Bbb C}^{q\times 2(p+q)},$ \ \
${\cal L}^{1,2}(z) \in {\Bbb C}^{2(p+q)\times p}.$

\section{Generalized inversion}

In the section we will present our main result.

Let us denote
\[
T_{{\cal R}_j}= \left(
\begin{array}{ccc}
{\cal R}_0^j & \ldots & 0 \\
\vdots &\ddots & \vdots \\
{\cal R}_m^j & \ldots & {\cal R}_0^j
\end{array}
\right), \ \
T_{{\cal L}_j}= \left(
\begin{array}{ccc}
{\cal L}_0^j & \ldots &  {\cal L}_{-n}^j \\
\vdots &\ddots & \vdots \\
0 & \ldots & {\cal L}_{0}^j
\end{array}
\right), \ \ j=1,2,
\]
where ${\cal R}_k^j ({\cal L}_k^j)$ are the coefficients of the polynomials ${\cal R}^j ({\cal L}^j)$.
We also put  $H_{{\cal R}_2}=JT_{{\cal R}_2},$ $H_{{\cal L}_1}=T_{{\cal L}_1}J$.

\begin{thm}
Generalized inverses of the $T+H$ and $T-H$ matrices are found by
the formulas:
\begin{equation}
\label{inv}
\left(T\pm H\right)^{\dagger}=\frac{1}{2}\left(T_{{\cal R}_1}\pm H_{{\cal R}_2}\right)\Pi
\left(T_{{\cal L}_2}\pm H_{{\cal L}_1}\right).
\end{equation}

If $T\pm H$ is invertible (one-sided invertible),
then $\left(T\pm H\right)^{\dagger}$
is its inverse (one-sided inverse) matrix.
\end{thm}

{\bf Proof.}

Let us construct a generalized inversion to $T_A\equiv T_0$ according to
formula (\ref{513}).

We are going to pass from block Toeplitz matrix $T_A$ to the mosaic matrix
\[
M_A=
\left(
\begin{array}{ccccccc}
b_n & \ldots& b_{n+m} &\vline& a_{n-m} & \ldots& a_{n} \\
b_{n-1} & \ldots &b_{n+m-1}&\vline& a_{n-m-1} & \ldots &a_{n-1}\\
\vdots & \ddots &\vdots & \vline&\vdots & \ddots &\vdots \\
b_{0} & \ldots & b_m& \vline&a_{-m} & \ldots & a_0\\
\hline
a_{0} &  \ldots& a_{-m}&\vline &b_m&  \ldots& b_0\\
a_{1} & \ldots &a_{-m+1}&\vline&b_{m+1}&\ldots &b_1\\
\vdots & \ddots &\vdots &\vline&\vdots & \ddots &\vdots \\
a_{n} & \ldots & a_{n-m}&\vline& b_{n+m}&\ldots &b_{n}
\end{array}
\right).
\]

At first,
according to the block structure of  $A_j$ (\ref{A}), we partition each
block column $X_j$ of the matrix $T_A$ into two block columns $X_j^1, X_j^2$ with sizes $2p(n+1)\times
q$:
\[
X_j=
\left(
\begin{array}{cc}
X_j^1 & X_j^2
\end{array}
\right).
\]

Then permute new block columns in $T_A$ and construct the matrix
\[
\left(
\begin{array}{cccccc}
X_1^1 &\ldots& X_m^1 & X_1^2 &\ldots & X_m^2
\end{array}
\right)=
\]
\[
\left(
\begin{array}{ccccccccc}
b_{n}&b_{n+1}&\ldots& b_{n+m} & \vline& a_{n-m} &a_{n-m+1}& \ldots & a_{n}\\
a_{0} &a_{-1}&\ldots& a_{-m} &\vline& b_{m} &b_{m-1}& \ldots & b_{0}\\
\hline
b_{n-1} &b_{n}&\ldots& b_{n+m-1} &\vline& a_{n-m-1} &a_{n-m}& \ldots & a_{n-1}\\
a_{1} &a_{0}&\ldots& a_{-m+1} & \vline& b_{m+1} &b_{m}& \ldots & b_{1}\\
\hline
\vdots & \vdots &\ddots &\vdots &\vline &\vdots &\vdots &\ddots &\vdots \\
\hline
b_{0} &b_{1}&\ldots& b_{m}& \vline & a_{-m} &a_{-m+1}& \ldots & a_{0}\\
a_{n} &a_{n-1}&\ldots& a_{n-m} &\vline& b_{n+m}  &b_{n+m-1}& \ldots & b_{n}
\end{array}
\right).
\]

This matrix is evident to be obtained by multiplying $T_A$ on a permutation matrix
$P_2$. Then we will do the analogous permutation with block rows in $T_AP_2$.
As a result, we will get the matrix $P_1T_AP_2$, where $P_1$ is a permutation
matrix. The matrix $P_1T_AP_2$ coincides with $M_A:$
\[M_A=P_1T_AP_2.
\]
Thus we have passed from the block Toeplitz matrix $T_A$ to the mosaic matrix
$M_A$.

Taking into account that for a permutation matrix $P$ the equality $P^{-1}~=~P^t$
holds, we get the generalized inversion for $M_A:$
\[
M_A^{\dagger}= P_2^tT_A^{\dagger}P_1^t.
\]

Let us specify the structure of factors in $M_A^{\dagger}= P_2^tT_A^{\dagger}P_1^t$.
The operations which $P_2$ has done with the block columns of $T_A$,
the matrix $P_2^t$ now will carry out with the block rows of the matrix
\[
\left(
\begin{array}{ccc}
{\cal R}_0 & \ldots & 0 \\
\vdots &\ddots & \vdots \\
{\cal R}_m & \ldots & {\cal R}_0
\end{array}
\right).
\]

Thus
\[
P_2^t\left(
\begin{array}{ccc}
{\cal R}_0 & \ldots & 0 \\
\vdots &\ddots & \vdots \\
{\cal R}_m & \ldots & {\cal R}_0
\end{array}
\right)=
\left(
\begin{array}{ccc}
{\cal R}_0^1 & \ldots & 0 \\
\vdots &\ddots & \vdots \\
{\cal R}_m^1 & \ldots & {\cal R}_0^1\\
{\cal R}_0^2 & \ldots & 0 \\
\vdots &\ddots & \vdots \\
{\cal R}_m^2 & \ldots & {\cal R}_0^2
\end{array}
\right)\equiv
\left(
\begin{array}{c}
T_{{\cal R}_1}\\
T_{{\cal R}_2}
\end{array}
\right),
\]
where ${\cal R}_j^{1,2}$ are the coefficients of the matrix polynomials
${\cal R}^{1,2}(z)$, presented in (\ref{div}).

Similarly, we have $\left(
\begin{array}{ccc}
{\cal L}_0 & \ldots &  {\cal L}_{-n} \\
\vdots &\ddots & \vdots \\
0 & \ldots & {\cal L}_{0}
\end{array}
\right)
P_1^t=$
\[=\left(
\begin{array}{cccccc}
{\cal L}_0^1 & \ldots &  {\cal L}_{-n}^1 & {\cal L}_0^2 & \ldots &  {\cal L}_{-n}^2\\
\vdots &\ddots & \vdots & \vdots &\ddots & \vdots \\
0 & \ldots & {\cal L}_{0}^1 & 0 & \ldots & {\cal L}_{0}^2
\end{array}
\right)\equiv\left(
\begin{array}{cc}
T_{{\cal L}_1}& T_{{\cal L}_2}
\end{array}\right).
\]

Then
\[
M_A^{\dagger}=\left(
\begin{array}{c}
T_{{\cal R}_1}\\
T_{{\cal R}_2}
\end{array}
\right)\Pi\left(
\begin{array}{cc}
T_{{\cal L}_1}& T_{{\cal L}_2}
\end{array}
\right).
\]

Let us apply now the well-known method \cite{Merchant}
of reducing the mosaic matrix $M_A$
to the block-diagonal matrix formed from the Toeplitz-plus-Hankel and Toeplitz-minus-Hankel
matrices:
\[
M_A=\frac{1}{2}\left(
\begin{array}{cc}
J & J \\
I & -I
\end{array}
\right)
\left(
\begin{array}{cc}
T+H & 0 \\
0 & T-H
\end{array}
\right)
\left(
\begin{array}{cc}
I & J \\
-I & J
\end{array}
\right).
\]

We obtain that the matrix
$G=\frac{1}{2}\left(
\begin{array}{cc}
I & J \\
-I & J
\end{array}
\right)M_A^{\dagger}
\left(
\begin{array}{cc}
J & J \\
I & -I
\end{array}
\right)=
$
\[=\frac{1}{2}\left(
\begin{array}{c}
T_{{\cal R}_1}+JT_{{\cal R}_2} \\
-T_{{\cal R}_1}+JT_{{\cal R}_2}
\end{array}
\right)\Pi
\left(
\begin{array}{cc}
T_{{\cal L}_1}J+T_{{\cal L}_2} & T_{{\cal L}_1}J-T_{{\cal L}_2}
\end{array}
\right)
\]
is the generalized inversion for the matrix $\left(
\begin{array}{cc}
T+H & 0 \\
0 & T-H
\end{array}
\right)$.

Let us present $G$ in the block form
\[
G=\left(
\begin{array}{cc}
G_{11} & G_{12} \\
G_{21} & G_{22}
\end{array}
\right),\]
with $G_{ij}\in {\Bbb C}^{(m+1)q\times (n+1)p}.$

It is easy to get that
\[
G_{11}=\frac{1}{2}\left(T_{{\cal R}_1}+ H_{{\cal R}_2}\right)\Pi
\left(T_{{\cal L}_2}+ H_{{\cal L}_1}\right),
\]
\[
G_{22}=\frac{1}{2}\left(T_{{\cal R}_1}- H_{{\cal R}_2}\right)\Pi
\left(T_{{\cal L}_2}- H_{{\cal L}_1}\right)
\]
are generalized inverses to $T+H$  $(T-H)$.

The theorem statement concerning the invertibility (one-sided invertibility)
is evident.

The theorem has been proved.

\halmos

Given $T\pm H$ matrices are block matrices with the sizes of their blocks
$p\times q.$
The factors in the inverse formulas (\ref{inv}) have blocks with sizes $q\times 2(p+q)$, $2(p+q)\times p$.
The compact form of the generalized inversion
is in many respects because of such factors sizes.
Sometimes it is convenient to have a formula for a generilized inversion
where factors have blocks with sizes $q\times q,q\times
p, p\times p$.

In order to obtain it we partition
${\cal R}(z)=\left(R_1(z) \ldots R_{2(p+q)}(z)\right)$ into blocks:
\[
{\cal R}(z)=\left(
\begin{array}{cccc}
{\cal R}_{11}& {\cal R}_{12}&{\cal R}_{13}& {\cal R}_{14}\\
{\cal R}_{21}&{\cal R}_{22}&{\cal R}_{23}&{\cal R}_{24}
\end{array}
\right).
\]

Here ${\cal R}_{ij}$ have the sizes $q\times q$ for $i,j=1,2,$ and
$q\times p$ for $i=1,2,$ $j=3,4.$

Then we will do the analogous partition with the matrix of the left
essential polynomials:
\[
{\cal L}(z)=\left(
\begin{array}{cc}
{\cal L}_{11}& {\cal L}_{12}\\
{\cal L}_{21}&{\cal L}_{22}\\
{\cal L}_{31}& {\cal L}_{32}\\
{\cal L}_{41}&{\cal L}_{42}
\end{array}
\right),
\]
here ${\cal L}_{ij}$ have the sizes $q\times p$ for $i,j=1,2,$ and $p\times p$ for
$i=3,4, j=1,2.$

Let us also partition
$D={\rm diag}[z^{\mu_1}\ldots z^{\mu_{2(p+q)}}]=\left(d_1\  d_2 \ d_3\  d_4\right),
$
where $d_{1,2}$ are diagonal matrices with the sizes $q\times q$ and $d_{3,4}$
are ones with the sizes $p\times p.$

We denote for $i,j=1,\ldots,4$
\[
T_{{\cal R}_{ij}}= \left(
\begin{array}{ccc}
{\cal R}_0^{ij} & \ldots & 0 \\
\vdots &\ddots & \vdots \\
{\cal R}_m^{ij} & \ldots & {\cal R}_0^j
\end{array}
\right), \ \
T_{{\cal L}_{ij}}= \left(
\begin{array}{ccc}
{\cal L}_0^{ij} & \ldots &  {\cal L}_{-n}^{ij} \\
\vdots &\ddots & \vdots \\
0 & \ldots & {\cal L}_{0}^{ij}
\end{array}
\right). \ \
\]

Then it is easy to see that
\[
\left(T\pm H\right)^{\dagger}=
\frac{1}{2}
\left[\sum_{j=1}^4 T_{{\cal R}_{1j}}\pi_jT_{{\cal L}_{j2}}+
\sum_{j=1}^4 H_{{\cal R}_{2j}}\pi_jH_{{\cal L}_{j1}}\pm
\right.
\]
\[\left.
\pm\left(
\sum_{j=1}^4 T_{{\cal R}_{1j}}\pi_jH_{{\cal L}_{j1}}+
\sum_{j=1}^4 H_{{\cal R}_{2j}}\pi_jT_{{\cal L}_{j2}}\right)\right],
\]
where we denote $T_{\cal L,R}J=H_{\cal L,R}$ and $\pi_j$ are the matrices
constructed by $d_j$ with the same manner as $\Pi$ by $d$.

\section{An example}

Let us demonstrate the theorem  with an example. We will find generalized inversions
of the following $T+H$ and $T-H$ matrices:
\[
T\pm H=\left(
\begin{array}{cccc}
1 & 0& -1& 1\\
1 & 1& 0 &-1\\
1 & 1 & 1 & 0\\
-1 & 1 &1 &1
\end{array}
\right)\pm \left(
\begin{array}{cccc}
1 & 0& -1& 1\\
0 & -1& 1 &0\\
-1 & 1 & 0 & 0\\
1 & 0 &0 &1
\end{array}
\right).
\]

Our calculations show that the indices of the sequence $A_{-3},\ldots,A_3$
are equal to $\mu_1=-1,\ \mu_2=0,\ \mu_3=0,\ \mu_4=1$. The matrix ${\cal R}(z)$
of the right essential polynomials and  the matrix ${\cal L}(z)$ of the conforming left
essential polynomials are:
\[
{\cal R}(z)=
\left(
\begin{array}{c}
{\cal R}^1(z) \\
{\cal R}^2(z)
\end{array}
\right), \ \ \
{\cal L}(z)=
\left(
\begin{array}{cc}
{\cal L}^1(z) & {\cal L}^2(z)
\end{array}
\right),
\]
where
\[
{\cal R}^1(z)=-\left(
\begin{array}{cccc}
1 &  2 & 2 &0
\end{array}\right)+
\left(
\begin{array}{cccc}
-4 &  5 & 3 &-1
\end{array}\right)z+\left(
\begin{array}{cccc}
2 &  1 & 0 &0
\end{array}\right)z^2-
\]
\[-\left(
\begin{array}{cccc}
11 &  4 & 2 &0
\end{array}\right)z^3+
\left(
\begin{array}{cccc}
0 &  4 & 1 &0
\end{array}\right)z^4+
\left(
\begin{array}{cccc}
0 &  0 & 0 &1
\end{array}\right)z^5,
\]
\[
{\cal R}^2(z)=\left(
\begin{array}{cccc}
11 &  -4 & -2 &2
\end{array}\right)-\left(
\begin{array}{cccc}
2 &  11 & 5 &-2
\end{array}\right)z+\left(
\begin{array}{cccc}
4 &  0 & 0 &0
\end{array}\right)z^2+
\]
\[+\left(
\begin{array}{cccc}
1 &  0 & 0 &0
\end{array}\right)z^3+
\left(
\begin{array}{cccc}
0 &  0& 1 &0
\end{array}\right)z^4,
\]
\[
{\cal L}^1(z)=\frac{1}{180}\left[\left(\begin{array}{c}
-4 \\ -12 \\ 4 \\ 0
\end{array}
\right)+\left(\begin{array}{c}
4 \\ 0 \\ 12 \\ 0
\end{array}
\right)z^{-1}+\left(\begin{array}{c}
-1 \\ 6 \\ -23 \\ 0
\end{array}
\right)z^{-2}-\right.
\]
\[\left.-\left(\begin{array}{c}
21 \\ -6 \\ 35 \\ 180
\end{array}
\right)z^{-3}+
\left(\begin{array}{c}
-18 \\ -84 \\ 76 \\ 0
\end{array}
\right)z^{-4}+\left(\begin{array}{c}
16 \\ 0 \\ 0 \\ 0
\end{array}
\right)z^{-5}
\right],
\]
\[
{\cal L}^2(z)=\frac{1}{180}\left[\left(\begin{array}{c}
25 \\ 30 \\ -25 \\ 180
\end{array}
\right)+\left(\begin{array}{c}
26 \\ 108 \\ -126 \\ 0
\end{array}
\right)z^{-1}+\left(\begin{array}{c}
-11 \\ 30 \\ -43 \\ 0
\end{array}
\right)z^{-2}+\right.
\]
\[\left.+\left(\begin{array}{c}
-6 \\ -24 \\ 32 \\ 0
\end{array}
\right)z^{-3}-
\left(\begin{array}{c}
10 \\ 24 \\ 26 \\ 0
\end{array}
\right)z^{-4}-\left(\begin{array}{c}
4 \\ 0 \\ 0 \\ 0
\end{array}
\right)z^{-5}
\right].
\]

We construct the matrix \[
\Pi=\left(
\begin{array}{cccc}
\Pi_0 & \Pi_{-1} & \Pi_{-2} & \Pi_{-3}\\
\Pi_1 & \Pi_0 & \Pi_{-1} & \Pi_{-2}\\
\Pi_2 & \Pi_1 & \Pi_0 & \Pi_{-1}\\
\Pi_3 & \Pi_2 & \Pi_1 & \Pi_0
\end{array}
\right),
\]
where
\[
\Pi_0=\left(
\begin{array}{cccc}
0&0&0&0\\
0&1&0&0\\
0&0&1&0\\
0&0&0&0
\end{array}
\right), \ \
\Pi_{-1}=\left(
\begin{array}{cccc}
0&0&0&0\\
0&0&0&0\\
0&0&0&0\\
0&0&0&1
\end{array}
\right),
\]
\[
\Pi_1=\left(
\begin{array}{cccc}
1&0&0&0\\
0&0&0&0\\
0&0&0&0\\
0&0&0&0
\end{array}
\right),
\]
and the matrices $\Pi_{\pm 2},\ \Pi_{\pm 3}$ are zero matrices.

Then we construct the matrices $T_{{\cal R}_1}$, $T_{{\cal R}_2}$ with
help of the coefficients of ${\cal R}^1,\ {\cal R}^2$ and matrices $T_{{\cal L}_1},$ $T_{{\cal L}_2}$
with help of the coefficients of ${\cal L}^1,\ {\cal L}^2$.

We use the formulas (\ref{inv}) for generalized inverses of the $T\pm H$ matrices
and have
\[
\left(T+H\right)^{\dagger}=\frac{1}{20}\left(
\begin{array}{cccc}
5& 10& 0& 0\\
2& -4& 12& -4\\
-4& 8& -4& 8\\
1& -2& -4& 8
\end{array}
\right),
\]
\[
\left(T-H\right)^{\dagger}=\frac{1}{180}\left(
\begin{array}{cccc}
-113& 2& 50& -32\\
88& 8& 20& 52\\
-134& -4& 80& 64\\
17& -158& 10 &16
\end{array}
\right).
\]

One may conclude that $\left(T+H\right)^{\dagger}$
is the inverse matrix for the $T+H$ matrix and $\left(T-H\right)^{\dagger}$
is the generalized inverse for $T-H.$

\end{document}